\documentclass[12pt,a4paper]{article}%
\usepackage{amsmath,amsfonts,eucal}%
\usepackage{amsmath}%
\usepackage{amsfonts}%
\usepackage{amssymb}%
\usepackage{graphicx} 
\usepackage{pstricks,textcomp,amssymb,amsmath,graphicx}

\begin{document}

\title{Solution of large linear systems
with embedded  network structure  for a non-homogeneous network
flow programming problem}
\author{L.A.~Pilipchuk, E.S.~Vecharynski}
\maketitle

\begin{abstract}
In the paper we consider the linear underdetermined system of a
special type. Systems of this type appear in non-homogeneous
network flow programming problems in the form of systems of
constraints and can be characterized as systems with a large
sparse submatrix representing the embedded network structure. We
develop a direct method for finding solutions of the system. The
algorithm is based on the theoretic-graph specificities for the
structure of the support and properties of the basis of a solution
space of a homogeneous system. One of the key steps is
decomposition of the system. A simple example is regarded at the
end of the paper.
\end{abstract}



\vspace{-0.2cm}\section*{\large \bf 1 Introduction}

\hspace*{\parindent}
        The work on this paper was motivated, mainly, by the analysis
of problems of non-homogeneous network flow optimization on large
data files [1]-[3], [5]-[7]. Our main goal was to develop an
effective (direct) method for solving large sparse systems of
linear equations with embedded network structure, which appear
naturally, e.g. as systems of constraints, in a broad class of
non-homogeneous network flow programming problems.

The 'network nature' of the regarded system allows keeping data in
the matrix-free form in the computer memory. The formulae, derived
within the paper, are written in the component (network) form to
provide clear approaches towards developing computational
algorithms using efficient data structures for graph
representation [1].

The general idea of the method is based on the following key
steps:

$\bullet$ {\it Distinguishing between the network part of the
system and the additional part.} The network part of the system
represents a network structure and corresponds to the network part
of the system of main constraints of a non-homogeneous network
flow programming problem [1], and is given, traditionally, by
balance equations, written for the nodes of a network. The
additional part of the system corresponds to the additional part
of the system of main constraints and can have a general form. We
start the solution by considering the network part of the system
only.

$\bullet$ {\it Introduction of the support of the network for a
system.} The term 'support of the network' {\it(also referred to
as network support, or support)} is borrowed from optimization
theory [2], [3] and is used here for further compatibility with
applications in problems of non-homogeneous network flow
programming. The actual meaning in this paper is -- a set of
indices of variables (or, in the network terms, - a set of arcs)
corresponding to columns, which form a basis minor of the matrix
of a system. We study the support for the network part of the
system, finding the correspondence between the columns of a basis
minor and a family of spanning trees.

$\bullet$ {\it  Construction of a general solution for the network
part of the system.} We compute a basis of a solution space of the
corresponding homogeneous system and interpret the basis vectors
as characteristic vectors, entailed by non-support arcs. A simple
approach for finding a partial solution of the (non-homogeneous)
system is provided.

$\bullet$ {\it  Decomposition of the system.}  We perform column
decomposition of the system by separating the variables according
to the sets - $U_{T},\,U_{C}$    and $U_{N}$, which consist of the
arcs of the support for the network part of the system, cyclic
arcs and non-support/non-cyclic arcs respectively; and, finally,
sequentially express the unknowns corresponding to the sets
$U_{C}$ and $U_{T}$  in terms of the independent variables
corresponding to the set  $U_{N}$.

\vspace{-0.2cm}\section*{\large \bf 1.1 General form of the
system}

\hspace*{\parindent}
 Let $S=(I,U)$  be a finite oriented connected network without
multiple arcs and loops, where $I$  is a set of nodes and $U$ is a
set of arcs defined on $I\times I(|I|<\infty,\,|U|<\infty)$ . Let
$K\,(|K|<\infty)$  be a set of different products (types of flow)
transported through the network $S$. For definiteness, we assume
the set $K=\{1,\ldots,|K|\}$. Let us denote a connected network
corresponding to a certain type of flow $k\in K$  with
${S}^{k}=(I^{k},{U}^{k}),I^{k}\subseteq
I,{U}^{k}=\{(i,j)^{k}:(i,j)\in\widetilde{U}^{k}\},\widetilde{U}^{k}
\subseteq U$ - a set of arcs of the network $S$ carrying the flow
of type $k$. Also, we define sets $K(i)=\{k\in K:i\in I^{k}\}$ and
$K(i,j)=\{k\in K:(i,j)^{k}\in {U}^{k}\}$ of types of flow
transported through a node $i\in I$ and an arc $(i,j)\in U$
respectively.

    Let us introduce a subset $U_{0}$  of the set $U$, and let
$K_{0}(i,j)\subseteq K(i,j),(i,j)\in U_{0}$  be an arbitrary
subset of $K(i,j)$ such that  $|K_{0}(i,j)|>1$.

    Finally, the initial network  $S=(I,U)$ may be considered as a union of
$|K|$ networks $S^{k}$, combined under additional constraints of a
general kind.

    Consider the following linear underdetermined system
\begin{equation}\label{eq1}
 \sum_{j\in I^{+}_{i}(U^{k})}x^{k}_{ij}-
\sum_{j\in I^{-}_{i}(U^{k})}x^{k}_{ji}=a^{k}_{i}, \quad i\in
I^{k},\,k\in K,
\end{equation}
\begin{equation}\label{eq2}
 \sum_{(i,j)\in U}
\sum_{k\in
K(i,j)}\lambda_{ij}^{kp}x^{k}_{ij}=\alpha_{p},\,p=\overline{1,q},
\end{equation}
\begin{equation}\label{eq3}
\sum_{k\in K_{0}(i,j)}x_{ij}^{k}=z_{ij},\quad(i,j)\in U_{0},
\end{equation}
where $I_{i}^{+}(U^{k})=\{j\in I^{k}: (i,j)^{k}\in U^{k}\}$ ,
$I_{i}^{-}(U^{k})=\{j\in I^{k}: (j,i)^{k}\in U^{k}\}$;
$a_{i}^{k},\lambda_{ij}^{kp},\alpha_{p},z_{ij}\in \mathbf{R}$ -
parameters of the system; $x=(x^{k}_{ij},(i,j)^{k}\in U^{k},k\in
K)$- vector of unknowns.


    The matrix of system (\ref{eq1}) - (\ref{eq3}) has the following block structure:
\begin{equation}\label{eq4}
A=\left[ \begin{array}{c}  M\\  Q\\  T\\
\end{array}
\right].
\end{equation}

Here $M$ is a sparse submatrix with a block-diagonal structure  of
size $\displaystyle\sum_{k\in K}|I^{k}|\times
\displaystyle\sum_{k\in K}|U^{k}|$ such that each block represents
a $|I^{k}|\times |U^{k}|$ incidence matrix of the network
$S^{k}=(I^{k},U^{k}),\,k\in K$, namely, \\* $M=M_{1}\bigoplus
M_{2}\bigoplus\cdots\bigoplus M_{|K|}$ , where
$M_{k},k=1,\ldots,|K|$ are blocks of matrix $M$; $Q$ is a $q\times
\displaystyle\sum_{k\in K}|U^{k}|$ submatrix (dense, in the
general case) with elements $\lambda_{ij}^{kp},(i,j)\in U$,$k\in
K(i,j)$ ,$p=\overline{1,q}$  ; $T$ is a
$|U_{0}|\times\displaystyle \sum_{k\in K}|U^{k}|$ submatrix
consisting of zeros and ones, where all the nonzero elements
appear in columns corresponding to arcs $(i,j)^{k},(i,j)\in
U_{0},k\in K_{0}(i,j)$. We assume that $\displaystyle\sum_{k\in
K}|I^{k}|+q+|U_{0}|<\displaystyle\sum_{k\in K}|U^{k}|$.

\vspace{-0.2cm}\section*{\large \bf 2 Network part of the system}

\hspace*{\parindent}
We start the solution of system (1) - (3) by
considering the network part of the system.

{\bf Definition 1}  {\it We call system (1) the network part of
the system (1)-(3). Systems (2) and (3) are called the additional
part of the system (1)-(3).}

Before we proceed, let us recall the following necessary and
sufficient condition of consistency for system (1) implied by
Kronecker-Capelli theorem: $$\sum_{i\in I^{k}}a^{k}_{i}=0, k\in
K.$$

{\bf Theorem 1. (Rank theorem)}.  The rank of the matrix of system
(\ref{eq1}) for the network  $S=(I,U)$ equals $\displaystyle
\sum_{k\in K}|I^{k}|-|K|$.

{\bf Proof.}  Since matrix  $M$ of the system (1) has the form \\
$M=M_{1}\bigoplus M_{2}\bigoplus\cdots\bigoplus M_{|K|}$, where
$M_{k}$ is a diagonal block of matrix  $M,\, k=1,\ldots,|K|$ and
 $rank\, M_{k}=|I^{k}|-1$ [1] then
$rank\, M=\displaystyle\sum^{|K|}_{k=1}rank\,
 M_{k}=\displaystyle
\sum_{k\in K}(|I^{k}|-1)=\displaystyle \sum_{k\in K}|I^{k}|-|K|$.
$\hfill \Box$

{\it Remark 1} We assume, without loss of generality, that the
rank of the system (\ref{eq1}) - (\ref{eq3}) \\ is $\displaystyle
\sum_{k\in K}|I^{k}|-|K|+q+|U_{0}|$ , where $q+|U_{0}|$  is a
number of equations in the additional part (2) - (3).

Since the matrix of system (\ref{eq1}) has the block-diagonal
structure, we split the solution of the system into $|K|$
solutions of (independent) systems , each of which corresponds to
a separate block, i.e. to a fixed $k\in K$ , and has the following
form:

\begin{equation}\label{eq5}
 \sum_{j\in I^{+}_{i}(U^{k})}x^{k}_{ij}-
\sum_{j\in I^{-}_{i}(U^{k})}x^{k}_{ji}=a^{k}_{i}, \quad i\in
I^{k}\quad.
\end{equation}

\vspace{-0.2cm}\section*{\large \bf 2.1 Support Criterion}

\hspace*{\parindent}
Let's define a support of the network
$S=(I,U)$ for system (\ref{eq1}).

{\bf Definition 2}  {\it The support of the network  $S=(I,U)$ for
system (\ref{eq1}) is a set of arcs  $U_{T}=\{U_{T}^{k}\subseteq
U^{k},k\in K\}$ such that the system
\begin{equation}\label{eq6}
 \sum_{j\in I^{+}_{i}(\hat{U}^{k})}x^{k}_{ij}-
\sum_{j\in I^{-}_{i}(\hat{U}^{k})}x^{k}_{ji}=0, \quad i\in
I^{k},\,k\in K
\end{equation}
has only a trivial solution for $\hat{U^{k}}=U^{k}_{T}$, but has a
non-trivial solution for $\hat{U^{k}}=U^{k}_{T},k\in K \setminus
k_{0};\,\hat{U^{k_{0}}}=U^{k_{0}}_{T}\bigcup
(i,j)^{k_{0}},\,(i,j)^{k_{0}}\not\in
U^{k_{0}}_{T},\,k_{0}\in K$.}\\[1mm]

{\bf Theorem 2. (Network Support Criterion).}  The set  \\
$U_{T}=\{U_{T}^{k},k\in K\}$ is a support of the network $S=(I,U)$
for system (\ref{eq1}) iff for each  $k\in K$ the set of arcs
$U^{k}_{T}$ is a spanning tree for the network
$S^{k}=(I^{k},U^{k})$.



{\bf Proof.}  Follows directly from the proof [3] for the case
when  $|K|=1$ and the block-diagonal structure of the matrix of
the system (\ref{eq1}). $\hfill \Box$

\vspace{-0.2cm}\section*{\large \bf 2.2 Basis of a solution space
of a homogeneous system. Characteristic vectors}

\hspace*{\parindent}
Before introducing the definition of a
characteristic vector, let's analyze the structure of a network
obtained by appending an arbitrary arc  $(\tau,\rho)^{k}\in
U^{k}\setminus U^{k}_{T}$, where $k\in K$  is fixed, to the
support $U_{T}$.

For a fixed $k\in K$  we consider a network
$\hat{S}^{k}=(I^{k},U^{k}_{T}\bigcup
(\tau,\rho)^{k}),(\tau,\rho)^{k}\in U^{k}\setminus U^{k}_{T}$,
where the set $U^{k}_{T}$ is a spanning tree of the network
$S^{k}$ . Appending an arc $(\tau,\rho)^{k}\in U^{k}\setminus
U^{k}_{T}$  to the tree entails a unique cycle. We denote this
cycle with $L^{k}_{\tau\rho}$. The set
$Z_{k}=\{L^{k}_{\tau\rho},(\tau,\rho)^{k}\in U^{k}\setminus
U^{k}_{T}\}$ is the fundamental set of cycles with respect to the
spanning of the network $S^{k}$ [1].

Let's consider a cycle  $L^{k}_{\tau\rho}$, entailed by an arc
$(\tau,\rho)^{k}\in U^{k}\setminus U^{k}_{T}$. We define the
detour direction within the cycle  $L^{k}_{\tau\rho}$
corresponding to the arc  $(\tau,\rho)^{k}$.\\[1mm]

{\bf Definition 3} {\it We call an arc  $(i,j)^{k}\in
L^{k}_{\tau\rho}$, where $k\in K$  is fixed, a forward arc of the
cycle $L^{k}_{\tau\rho}$, if the direction of  the arc $(i,j)^{k}$
is the same as the direction of  the arc $(\tau,\rho)^{k}$  within
the cycle  $L^{k}_{\tau\rho}$. Similarly, we call an arc
$(i,j)^{k}\in L^{k}_{\tau\rho}$, where $k\in K$  is fixed, a
backward arc of the cycle  $L^{k}_{\tau\rho}$, if the direction of
the arc  $(i,j)^{k}$ is opposite to the direction of the arc
$(\tau,\rho)^{k}$ within the cycle $L^{k}_{\tau\rho}$.}

    We denote the sign of an arc $(i,j)^{k}$   within a cycle
$L^{k}_{\tau\rho}$  by $sign(i,j)^{L^{k}_{\tau\rho}}$ ,

\begin{equation}\label{eq7}
sign(i,j)^{L^{k}_{\tau\rho}}=\left\{ \begin{array}{c}
1,\,\, (i,j)^{k}\in L^{k+}_{\tau\rho}\\
-1,(i,j)^{k}\in L^{k-}_{\tau\rho},\\
0,\,(i,j)^{k}\not\in L^{k}_{\tau\rho}\\
\end{array}
\right.
\end{equation}
where  $L^{k_{+}}_{\tau\rho}$ and $L^{k_{-}}_{\tau\rho}$  are the
sets of forward and backward arcs of the cycle $L^{k}_{\tau\rho}$
with a direction corresponding to the arc  $(\tau,\rho)^{k}$.

Let us give a constructive definition of a characteristic vector,
entailed by an arc.

{\bf Definition 4} {\it Characteristic vector, entailed by an arc
$(\tau,\rho)^{k}\in U^{k}\setminus U^{k}_{T}$  with respect to the
spanning tree $U^{k}_{T}$ , is a vector
$\delta^{k}(\tau,\rho)=(\delta^{k}_{ij}(\tau,\rho),(i,j)^{k}\in
U^{k})$, where $k\in K$ is fixed, constructed according to the
following rules:

$\bullet$  Add an arc  $(\tau,\rho)^{k}\in U^{k}\setminus
U^{k}_{T},$ to the set $U^{k}_{T},k\in K$ , which is a spanning
tree for the network $S^{k}=(I^{k},U^{k})$; and thus create a
unique cycle $L^{k}_{\tau\rho}$.

$\bullet$   Let the arc $(\tau,\rho)^{k}$  set the detour
direction within the cycle $L^{k}_{\tau\rho}$ and
$\delta^{k}_{\tau\rho}(\tau,\rho)=1$.

$\bullet$ For cycle's forward arcs, let
$\delta^{k}_{ij}(\tau,\rho)=1$.

$\bullet$   For cycle's backward arcs, let
$\delta^{k}_{ij}(\tau,\rho)=-1$.

$\bullet$   Let $\delta^{k}_{ij}(\tau,\rho)=0,(i,j)^{k}\in
U^{k}\backslash L^{k}_{\tau\rho} $ .}

For briefness, further in this paper, we will call a
characteristic vector  $\delta^{k}(\tau,\rho)$, entailed by an arc
$(\tau,\rho)^{k}$,  with respect to the spanning tree $U_{T}^{k}$,
a characteristic vector  $\delta^{k}(\tau,\rho)$, entailed by an
arc  $(\tau,\rho)^{k}$, or, simply, a characteristic vector
$\delta^{k}(\tau,\rho)$.

The next two lemmas state the essential properties of
characteristic vectors.\\[1mm]

{\bf Lemma 1} A characteristic vector $\delta^{k}(\tau,\rho)$,
entailed by an arc $(\tau,\rho)^{k} \in U^{k}\setminus U^{k}_{T}$,
where $k\in K$  is fixed, is a solution of the homogeneous linear
system (\ref{eq8})

\begin{equation}\label{eq8}
 \sum_{j\in I^{+}_{i}(U^{k})}x^{k}_{ij}-
\sum_{j\in I^{-}_{i}(U^{k})}x^{k}_{ji}=0, \quad i\in I^{k}.
\end{equation}

{\bf Proof.} Let a support $U_{T}=\{U_{T}^{k},k\in K\}$  be
defined. For a fixed $k\in K$ we consider the set $U_{T}^{k}$
which is, according to Theorem 2, a spanning tree for the network
$S^{k}$, and let $L^{k}_{\tau\rho}$ be the unique cycle of the
network $\hat{S}^{k}=(I^{k},U^{k}_{T}\bigcup(\tau,\rho)^{k})$,
which appears after appending the arc $(\tau,\rho)^{k}\in
U^{k}\setminus U^{k}_{T}$  to the set $U^{k}_{T}$.

Consider the vector $x^{k}=(x^{k}_{ij},(i,j)^{k}\in U^{k})$ of
unknowns in system (\ref{eq8}).

Let's let $x^{k}_{ij}=0,(i,j)^{k}\in U^{k}\backslash
L^{k}_{\tau\rho}$. Thus, the system (\ref{eq8}) can be reduced to

\begin{equation}\label{eq9}
\displaystyle \sum_{j\in I^{+}_{i}(L^{k}_{\tau\rho})}x^{k}_{ij}-
\displaystyle \sum_{j\in I^{-}_{i}(L^{k}_{\tau\rho})}x^{k}_{ji}=0,
\quad i\in I(L^{k}_{\tau\rho}),
\end{equation}
where $I(L^{k}_{\tau\rho})$  denotes all nodes in cycle
$L^{k}_{\tau\rho}$.

Letting   $x^{k}_{\tau\rho}=1$, from the reduced system (9), we
can easily define the values of the remaining unknowns
$x^{k}_{ij},(i,j)^{k}\in
L^{k}_{\tau\rho}\setminus(\tau,\rho)^{k}$:

$$
x^{k}_{ij}=sign(i,j)^{L^{k}_{\tau\rho}},(i,j)^{k}\in
L^{k}_{\tau\rho}\setminus(\tau,\rho)^{k}.
$$

Algorithmically, after letting  $x^{k}_{\tau\rho}=1$, we pass from
node $\tau$  to node $\rho$ along the cycle  $L^{k}_{\tau\rho}$,
consecutively setting the unknowns  $x^{k}_{ij},(i,j)^{k}\in
L^{k}_{\tau\rho}\setminus(\tau,\rho)^{k}$,   to the values of
signs of the corresponding arcs within the cycle
$L^{k}_{\tau\rho}$.

Note, the constructed solution vector $x^{k}$ satisfies all the
rules of Definition 4 of a characteristic vector, entailed by an
arc $(\tau,\rho)^{k}$, and hence $\delta^{k}(\tau,\rho)=x^{k}$ is
a solution of the homogeneous linear system (\ref{eq8}).
$\hfill\Box$
\\[1mm]

{\bf Lemma 2.} The set $\{\delta^{k}(\tau,\rho),(\tau,\rho)^{k}\in
U^{k}\setminus U^{k}_{T}\}$ of characteristic vectors, where $k\in
K$ is fixed, forms the basis of a solution space for the
homogeneous system (\ref{eq8}).

{\bf Proof.} According to Lemma 1, each characteristic vector
satisfies the homogeneous system (\ref{eq8}).

By Theorem 2, for a fixed  $k\in K$, the set  $U^{k}_{T}$ is a
spanning tree for the network  $S^{k}=(I^{k},U^{k})$, hence
$|U^{k}_{T}|=|I^{k}|-1$. Thus, the number of characteristic
vectors in the set  $\{\delta^{k}(\tau,\rho),(\tau,\rho)^{k}\in
U^{k}\setminus U^{k}_{T}\}$  equals $|U^{k}\setminus
U^{k}_{T}|=|U^{k}|-|I^{k}|+1$.

Now it suffices to show that all the vectors in the set are
linearly independent.

Each characteristic vector  $\delta^{k}(\tau,\rho)$, entailed by
some arc  $(\tau,\rho)^{k}\in U^{k}\setminus U^{k}_{T}$, always
has one and only one component, corresponding to the set
$U^{k}\setminus U^{k}_{T}$, that is equal to 1. It corresponds to
the arc  $(\tau,\rho)^{k}\in U^{k}\setminus U^{k}_{T}$ that has
entailed this vector. All the other components, which correspond
to arcs  $U^{k}\setminus L^{k}_{\tau\rho}$, are equal to $0$. This
fact implies that any two characteristic vectors, entailed by
different arcs, are linearly independent.$\hfill\Box$
\\[1mm]

{\bf Theorem 3.} The general solution of system (\ref{eq5}), for a
fixed $k\in K $, can be represented using the following form:

\begin{gather}\label{eq10}
x^{k}_{ij}= \sum_{(\tau,\rho)^{k}\in U^{k}\backslash
U^{k}_{T}}x^{k}_{\tau\rho}sign(i,j)^{L^{k}_{\tau\rho}}
+\left (\tilde{x}^{k}_{ij}-\sum_{(\tau,\rho)^{k}\in
U^{k}\backslash U^{k}_{T}}
\tilde{x}^{k}_{\tau\rho}sign(i,j)^{L^{k}_{\tau\rho}}\right),\\
(i,j)^{k}\in U^{k}_{T},\quad x^{k}_{\tau\rho}\in \mathbf{R},\quad
(\tau,\rho)^{k}\in U^{k}\setminus U^{k}_{T}, \nonumber
\end{gather}
where $\tilde{x}^{k}=(\tilde{x}^{k}_{ij},(i,j)^{k}\in
U^{k})$ is any partial solution of the (non-homogeneous) system
(\ref{eq5});
 $x^{k}_{\tau\rho}$ are independent variables corresponding to
arcs $(\tau,\rho)^{k}\in U^{k}\setminus U^{k}_{T}$.

{\bf Proof.} Let $x^{k}=(x^{k}_{ij},(i,j)^{k}\in U^{k})$ be a
general solution, and \\
$\tilde{x}^{k}=(\tilde{x}^{k}_{ij},(i,j)^{k}\in U^{k})$ - a
partial solution, of the system (\ref{eq5}).  Since, by Lemma 2,
the set $\{\delta^{k}(\tau,\rho),(\tau,\rho)^{k}\in U^{k}\setminus
U^{k}_{T}\}$ of characteristic vectors forms the basis of a
solution space for the homogeneous system (\ref{eq8}), we can
write the expression for $x^{k}$  in the following vector form:
\begin{equation}\label{eq11}
x^{k}=\sum_{(\tau,\rho)^{k}\in U^{k}\backslash U^{k}_{T}}
\alpha^{k}_{\tau\rho}\delta^{k}(\tau,\rho) +\tilde{x}^{k},
\end{equation}
as a sum of a general solution of the homogeneous system
(\ref{eq8}) and a partial solution of the non-homogeneous system
(\ref{eq5}); $\alpha^{k}_{\tau\rho}\in \mathbf{R}$  are
coefficients of the linear combination of characteristic vectors
in (\ref{eq11}).

Rewriting (\ref{eq11}) in the component form we obtain:

\begin{equation}\label{eq12}
x^{k}_{ij}=\sum_{(\tau,\rho)^{k}\in U^{k}\backslash U^{k}_{T}}
\alpha^{k}_{\tau\rho}\delta^{k}_{ij}(\tau,\rho)
+\tilde{x}^{k}_{ij},\quad (i,j)^{k}\in U^{k}_{T};
\end{equation}
\begin{equation}\label{eq13}
x^{k}_{\tau\rho}=\alpha^{k}_{\tau\rho}
+\tilde{x}^{k}_{\tau\rho},\quad (\tau,\rho)^{k}\in U^{k}\backslash
U^{k}_{T}.
\end{equation}

From equations (\ref{eq13}) we find
$\alpha^{k}_{\tau\rho}=x^{k}_{\tau\rho}
-\tilde{x}^{k}_{\tau\rho},\quad (\tau,\rho)^{k}\in U^{k}\backslash
U^{k}_{T}$   and substitute into (\ref{eq12}). Finally, rewriting
components of characteristic vectors according to (\ref{eq7}), we
obtain the expression (10) for the general solution of the system
(\ref{eq5}). $\hfill\Box$

{\it Remark 2} In practice, for construction of a partial
solution\\
 $\tilde{x}^{k}=(\tilde{x}^{k}_{ij},(i,j)^{k}\in U^{k})$
of the system (\ref{eq5}), we a priori assume
$\tilde{x}^{k}_{\tau\rho}=0 $, $(\tau,\rho)^{k}\in U^{k}\backslash
U^{k}_{T}$ and solve the system
$$\displaystyle \sum_{j\in I^{+}_{i}(U^{k}_{T})}\tilde{x}^{k}_{ij}-
\displaystyle \sum_{j\in
I^{-}_{i}(U^{k}_{T})}\tilde{x}^{k}_{ji}=a^{k}_{i}, \quad i\in
I^{k}.$$

Thus, formula (\ref{eq10}) gets to the form:
\begin{equation}\label{eq14}
x^{k}_{ij}=\sum_{(\tau,\rho)^{k}\in U^{k}\backslash U^{k}_{T}}
x^{k}_{\tau\rho}sign(i,j)^{L^{k}_{\tau\rho}}
+\tilde{x}^{k}_{ij}\,,(i,j)^{k}\in U^{k}_{T},
\end{equation}
$$
 x^{k}_{\tau\rho}\in \mathbf{R},\quad
(\tau,\rho)^{k}\in U^{k}\backslash U^{k}_{T}.
$$

Further, we will use the formula (\ref{eq14}).

\vspace{-0.2cm}\section*{\large \bf 3 Decomposition of the system}

\hspace*{\parindent} Let $U_{T}=\{U_{T}^{k},k\in K\}$   be a
support of the network  $S$ for the system (\ref{eq1}). We define
a set $U_{C}=\{U_{C}^{k}\subseteq U^{k}\backslash U_{T}^{k},k\in
K\},|U_{C}|=q+|U_{0}|$ of cyclic arcs by selecting $q+|U_{0}|$
arbitrary arcs from the sets $U^{k}\backslash U_{T}^{k},k\in K$.
We denote $U_{N}=\{U_{N}^{k},k\in
K\},U_{N}^{k}=U^{k}\backslash(U_{T}^{k}\bigcup U_{C}^{k}),k\in K $
- the set of remaining arcs, which were not included neither to
the support  $U_{T}$, nor to the set of cyclic arcs $U_{C}$.

Let's substitute the general solution (\ref{eq14}) of the system
(\ref{eq5}), for each \\ $k\in K$,  into (\ref{eq2}):

$$\sum_{(i,j)\in U} \sum_{k\in
K(i,j)}\lambda_{ij}^{kp}x^{k}_{ij}=\sum_{k\in K}\sum_{(i,j)^{k}\in
U^{k}} \lambda_{ij}^{kp}x^{k}_{ij}=$$
$$= \sum_{k\in
K}\sum_{(i,j)^{k}\in U^{k}_{T}} \lambda_{ij}^{kp}x^{k}_{ij}+
\sum_{k\in K}\sum_{(\tau,\rho)^{k}\in U^{k}\backslash
U^{k}_{T}}\lambda_{\tau\rho}^{kp}x^{k}_{\tau\rho}=
$$
$$
=\sum_{k\in
K}\sum_{(i,j)^{k}\in U^{k}_{T}} \lambda_{ij}^{kp}\left[
\sum_{(\tau,\rho)^{k}\in U^{k}\backslash U^{k}_{T}}
x^{k}_{\tau\rho}sign(i,j)^{L^{k}_{\tau\rho}}
+\tilde{x}^{k}_{ij}\right]+
$$
\begin{equation}\label{eq15}
 +\sum_{k\in K}\sum_{(\tau,\rho)^{k}\in U^{k}\backslash
U^{k}_{T}}\lambda_{\tau\rho}^{kp}x^{k}_{\tau\rho}=\alpha_{p},\quad
p=\overline{1,q}
\end{equation}

    We change the summing order in (\ref{eq15}):
$$
\sum_{k\in K}\sum_{(\tau,\rho)^{k}\in U^{k}\backslash U^{k}_{T}}
x^{k}_{\tau\rho}\sum_{(i,j)^{k}\in
U^{k}_{T}}\lambda_{ij}^{kp}sign(i,j)^{L^{k}_{\tau\rho}}
+\sum_{k\in K}\sum_{(i,j)^{k}\in U^{k}_{T}}\lambda_{ij}^{kp}
\tilde{x}^{k}_{ij}+
$$
\begin{equation}\label{eq16}
+\sum_{k\in K}\sum_{(\tau,\rho)^{k}\in U^{k}\backslash U^{k}_{T}}
\lambda_{\tau\rho}^{kp}x^{k}_{\tau\rho}=\alpha_{p},\quad
p=\overline{1,q}\,.
\end{equation}

In equations (\ref{eq16}) we group the variables, corresponding to
the sets  $U^{k}\backslash U_{T}^{k}$, $ k\in K$:
$$
\sum_{k\in K}\sum_{(\tau,\rho)^{k}\in U^{k}\backslash U^{k}_{T}}
x^{k}_{\tau\rho}\left[\lambda_{\tau\rho}^{kp}+\sum_{(i,j)^{k}\in
U^{k}_{T}}\lambda_{ij}^{kp}sign(i,j)^{L^{k}_{\tau\rho}}\right]=
$$
\begin{equation}\label{eq17}
= \alpha_{p}-\sum_{k\in K}\sum_{(i,j)^{k}\in
U^{k}_{T}}\lambda_{ij}^{kp}
\tilde{x}^{k}_{ij}\,,\,p=\overline{1,q}\,.
\end{equation}

{\bf Definition 5} {\it We call the number
\begin{equation}\label{eq18}
R_{p}(L^{k}_{\tau\rho})=\displaystyle \sum_{(i,j)^{k}\in
L^{k}_{\tau\rho}}\lambda_{ij}^{kp}sign(i,j)^{L^{k}_{\tau\rho}}
\end{equation}
the determinant of the cycle  $L^{k}_{\tau\rho}$, entailed by an
arc  $(\tau,\rho)^{k}\in U^{k}\backslash U_{T}^{k}$, with respect
to the equation with the number
$p$  of the system (\ref{eq2}).}\\[1mm]

Let's denote

\begin{equation}\label{eq19}
A^{p}=\alpha_{p}-\sum_{k\in K}\sum_{(i,j)^{k}\in
U^{k}_{T}}\lambda_{ij}^{kp}\,
\tilde{x}^{k}_{ij}\,,\,p=\overline{1,q}\,.
\end{equation}

The equations (\ref{eq17}), according to formulae (\ref{eq18}),
(\ref{eq19}), get to the form:
\begin{equation}\label{eq20}
\sum_{k\in K}\sum_{(\tau,\rho)^{k}\in U^{k}\backslash U^{k}_{T}}
R_{p}(L^{k}_{\tau\rho})x^{k}_{\tau\rho}=
 A^{p},\quad p=\overline{1,q}\,.
\end{equation}

In (\ref{eq20}) we group the variables, corresponding to the sets
$U^{k}_{C},\,k\in K$:
\begin{equation}\label{eq21}
\sum_{k\in K}\sum_{(\tau,\rho)^{k}\in U^{k}_{C}}
R_{p}(L^{k}_{\tau\rho})x^{k}_{\tau\rho}=
 A^{p}-\sum_{k\in K}\sum_{(\tau,\rho)^{k}\in U^{k}_{N}}
R_{p}(L^{k}_{\tau\rho})x^{k}_{\tau\rho} ,p=\overline{1,q}\,.
\end{equation}

Now, we apply the similar considerations to the system
(\ref{eq3}). Note, that (\ref{eq3}) can be regarded as a
particular case of the system (\ref{eq2}) with $\lambda_{ij}^{kp}$
equal to 0 or 1.


Let us substitute the general solution (\ref{eq14}) of the system
(\ref{eq5}), for each $k\in K$ , into (\ref{eq3}):
$$
\sum_{ k\in K_{0}(i,j)} x^{k}_{ij} = \sum_{\scriptstyle k\in
K_{0}(i,j),\atop \scriptstyle(i,j)^{k}\in U^{k}_{T}} x^{k}_{ij} +
\sum_{\scriptstyle k\in K_{0}(i,j),\atop \scriptstyle(i,j)^{k}\in
U^{k}\backslash U^{k}_{T}} x^{k}_{ij} =
$$
$$
=\sum_{\scriptstyle k\in K_{0}(i,j),\atop \scriptstyle(i,j)^{k}\in
U^{k}_{T}} \left[ \sum_{(\tau,\rho)^{k}\in U^{k}\backslash
U^{k}_{T}}
x^{k}_{\tau\rho}sign(i,j)^{L^{k}_{\tau\rho}}+\tilde{x}^{k}_{ij}\right]+
$$
\begin{equation}\label{eq22}
+\sum_{\scriptstyle k\in K_{0}(i,j),\atop \scriptstyle(i,j)^{k}\in
U^{k}\backslash U^{k}_{T}} x^{k}_{ij} =z_{ij},\quad (i,j)\in
U_{0}.
\end{equation}

Now, after changing the summing order and grouping the variables,
corresponding to the sets $U^{k}\backslash U^{k}_{T},\,k\in
K_{0}(i,j),\,(i,j)\in U_{0}$  in (\ref{eq22}), we obtain

\begin{equation}\label{eq23}
\sum_{\scriptstyle k\in K_{0}(i,j),\atop
\scriptstyle(\tau,\rho)^{k}\in U^{k}\backslash U^{k}_{T}}
x^{k}_{\tau\rho} sign(i,j)^{L^{k}_{\tau\rho}}= z_{ij}-
\sum_{\scriptstyle k\in K_{0}(i,j),\atop \scriptstyle(i,j)^{k}\in
U^{k}_{T}}\tilde{x}^{k}_{ij}, \quad (i,j)\in U_{0}\,.
\end{equation}

On this step let us introduce the following notation:

\begin{equation}\label{eq24}
\delta_{ij}(L^{k}_{\tau\rho})=\left\{ \begin{array}{lc}
sign(i,j)^{L^{k}_{\tau\rho}},k\in K_{0}(i,j)&  \\
&,(i,j)\in U_{0},(\tau,\rho)^{k}\in U^{k}\backslash U^{k}_{T},k\in K.\\
0,k \not \in K_{0}(i,j)&
\end{array}
\right.
\end{equation}

Thus, equations (\ref{eq23}) get to the form
\begin{equation}\label{eq25}
 \sum_{k\in K}\sum_{(\tau,\rho)^{k}\in
U^{k}\backslash U^{k}_{T}}\delta_{ij}(L^{k}_{\tau\rho})
x^{k}_{\tau\rho}=A_{ij}, \quad (i,j)\in U_{0},
\end{equation}
where
\begin{equation}\label{eq26}
A_{ij}=z_{ij}- \sum_{\scriptstyle k\in K_{0}(i,j), \atop
\scriptstyle (i,j)^{k}\in U^{k}_{T}}\tilde{x}^{k}_{ij},
\quad(i,j)\in U_{0}.
\end{equation}

In (\ref{eq25}) we group the variables, corresponding to the sets
$U^{k}_{C},k\in K$:
$$
 \sum_{k\in K }\sum_{(\tau,\rho)^{k}\in
U^{k}_{C}}\delta_{ij}(L^{k}_{\tau\rho}) x^{k}_{\tau\rho}=
$$
\begin{equation}\label{eq27}
=A_{ij}-\sum_{k\in K }\sum_{(\tau,\rho)^{k}\in
 U^{k}_{N}}\delta_{ij}(L^{k}_{\tau\rho})
x^{k}_{\tau\rho}
 , \quad (i,j)\in U_{0}.
\end{equation}

Finally, let us rewrite equations (\ref{eq21}) and (\ref{eq27}) in
the matrix form. For this purpose, we introduce arbitrary
numberings of arcs within the sets $U_{0}$  and  $U_{C}$. Thus,
$\xi=\xi(i,j)$ is a number of an arc  $(i,j)\in U_{0}, \,\xi \in
\{1,2,\ldots,|U_{0}|\}$; and $t=t(\tau,\rho)^{k}$ is a number of a
cyclic arc $(\tau,\rho)^{k}\in U_{C}^{k},k\in K ,t\in
\{1,2,\ldots,|U_{C}|\}$. In other words, we number the equations
of the system (\ref{eq3}), or (\ref{eq27}), and the variables,
corresponding to the set $U_{C}$. Note, the numbering of cyclic
arcs is equivalent to the numbering of the set
$\{L^{k}_{\tau\rho},(\tau,\rho)^{k}\in U_{C}^{k},k\in K\}$  of
cycles, entailed by arcs $(\tau,\rho)^{k}\in U_{C}^{k}$, with
respect to spanning trees $U_{T}^{k}$ of the networks  $S^{k}$.

Now equations (\ref{eq21}) and (\ref{eq27}) can be regarded as
following:
\begin{equation}\label{eq28}
Dx_{C}=\beta,
\end{equation}
where $D=\left({D_{1}\atop D_{2}}\right)$,
$D_{1}=(R_{p}(L^{k}_{\tau\rho}),p=\overline{1,q},t(\tau,\rho)^{k}=
\overline{1,|U_{C}|})$ - submatrix of the size $q\times|U_{C}|$,
$D_{2}=(\delta_{ij}(L^{k}_{\tau\rho}),\xi(i,j)=\overline{1,|U_{0}|},
t(\tau,\rho)^{k}=\overline{1,|U_{C}|})$ - submatrix of the size
$|U_{0}|\times|U_{C}|$,
$x_{C}=(x^{k}_{\tau\rho},(\tau,\rho)^{k}\in U^{k}_{C},k\in K)$
 - vector of unknowns with components ordered according to the
numbering  $t=t(\tau,\rho)^{k}$.

The right-hand side of (\ref{eq28}) has the form:

\begin{equation}\label{eq29}
\beta=\left( \begin{array}{cc}
\beta_{p},& p=\overline{1,q} \\
\beta_{q+\xi(i,j)},& (i,j)\in U_{0} \\
\end{array}
\right),
\end{equation}
where  $\beta_{p}=A^{p}-\displaystyle \sum_{k\in
K}\sum_{(\tau,\rho)^{k}\in U^{k}_{N}}R_{p}(L^{k}_{\tau\rho})
x^{k}_{\tau\rho},p=\overline{1,q}$,
$$
\beta_{q+\xi(i,j)}=A_{ij}-\displaystyle \sum_{k\in
K}\sum_{(\tau,\rho)^{k}\in U^{k}_{N}}\delta_{ij}(L^{k}_{\tau\rho})
x^{k}_{\tau\rho}
 , \quad (i,j)\in U_{0} .
$$

From (\ref{eq28}), in case of non-singularity of the matrix  $D$,
we find the unknown variables  $x_{C}$, corresponding to the set
$U_{C}$ of cyclic arcs:
\begin{equation}\label{eq30}
x_{C}=D^{-1}\beta\,.
\end{equation}

{\it Remark 3.} Generally, because of an arbitrary selection of
arcs for the set  $U_{C}=\{U_{C}^{k},k\in K\}$, non-singularity of
the matrix $D$ is not guaranteed. In the case when det  $D=0$ one
should re-select arcs into the set $U_{C}$ and re-compute
$D,\beta$ for the system (\ref{eq28}).

Let  $D^{-1}=(\nu_{l,s};l,s=\overline{1,|U_{C}|})$. We rewrite
(\ref{eq30}) in the component form:
$$
x^{k}_{\tau\rho}=\sum_{p=1}^{q}\nu_{t,p}\beta_{p}+\sum_{(i,j)\in
U_{0}}\nu_{t,q+\xi(i,j)}\beta_{q+\xi(i,j)},\,t=t(\tau,\rho)^{k},
\,(\tau,\rho)^{k}\in U^{k}_{C},k\in K.
$$

Thus, we have determined all the unknown variables
$x^{k}=(x^{k}_{ij},(i,j)^{k}\in U^{k},k\in K)$  of the system
(\ref{eq1}) - (\ref{eq3}):
$$
x^{k}_{\tau\rho}=\sum_{p=1}^{q}\nu_{t,p}\beta_{p}+\sum_{(i,j)\in
U_{0}}\nu_{t,q+\xi(i,j)}\beta_{q+\xi(i,j)},
$$
\begin{equation}\label{eq31}
t=t(\tau,\rho)^{k}, \,(\tau,\rho)^{k}\in U^{k}_{C},\,k\in K,
\end{equation}
\begin{equation}\label{eq32}
x^{k}_{ij}=\sum_{(\tau,\rho)^{k}\in U^{k}_{N}}
x^{k}_{\tau\rho}sign(i,j)^{L^{k}_{\tau\rho}}+\psi^{k}_{ij}
+\tilde{x}^{k}_{ij},\, (i,j)^{k}\in U^{k}_{T},\,k\in K,
\end{equation}
$$x^{k}_{\tau\rho}\in
\mathbf{R}\,,\, (\tau,\rho)^{k}\in U^{k}_{N},$$ where
$\psi^{k}_{ij}=\displaystyle \sum_{(\tau,\rho)^{k}\in U^{k}_{C}}
x^{k}_{\tau\rho}sign(i,j)^{L^{k}_{\tau\rho}}$.

Note, the components of the vector
$\tilde{x}^{k}=(\tilde{x}^{k}_{ij},(i,j)^{k}\in U^{k})$  of a
partial solution of the system (\ref{eq5}) are constructed
according to the rules in the {\it Remark 2}.

Before we start with a simple example, let us briefly discuss the
most important, in our opinion, aspects of the method. Although
the strict estimate of complexity was left beyond the scope of the
paper, one can notice that the described approach, if implemented
on proper data structures, leads to efficient algorithm: the
reasonable part of computations is done on small subsets of arcs,
e.g. on 'isolated' cycles - (\ref{eq7}), (\ref{eq18}), or spanning
trees - (\ref{eq19}), (\ref{eq26}). The use of the embedded
network structure allows performing decomposition of the system
and, finally, inverting the matrix $D$ (\ref{eq28}) of a size much
smaller than that of the initial system (\ref{eq1})-(\ref{eq3}).
Moreover, the fact that the same results were obtained for each
type of flow $k\in K$, e.g. Theorem 2, Lemmas 1 and 2, formulae
(\ref{eq14}), makes the method ready for implementation in
parallel environment.

However, the power of the approach is appreciated in the context
of large problems of non-homogeneous network flow programming with
(\ref{eq1})-(\ref{eq3}) being the system of main constraints,
where the presented ideas provide the uniform technique for
computing essential quantities: increment of an objective
function, feasible directions, pseudo-flow, etc.

Currently the authors work on the application of the obtained
results for derivation of an optimality criterion for a broad
class of non-homogeneous network flow programming problems.

\vspace{-0.2cm}\section*{\large \bf 4  Example}

\begin{figure}[h]
    \centering
        \includegraphics[width=300pt,height=300pt]{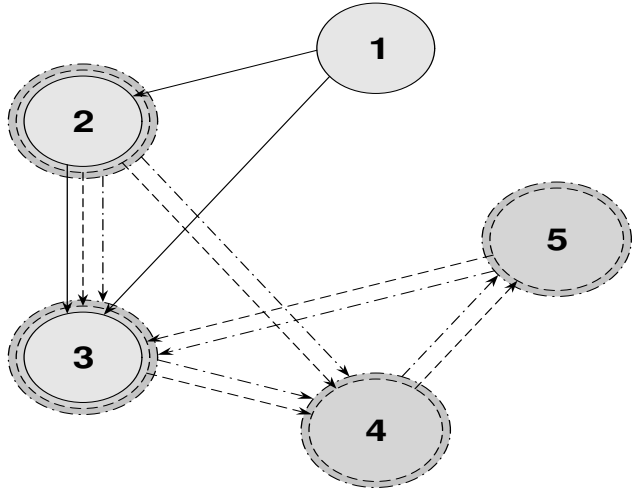}
    \caption{Union of networks $S^{k}=(I^{k},U^k), k\in K=\{1,2,3\}$}
    \label{graph}
\end{figure}

\hspace*{\parindent} Let us consider the example (1a) - (3a) of
the problem (\ref{eq1}) - (\ref{eq3}) for the network $S=(I,U), I
= \{1,2,3,4,5\}, U=\{(1,2),(1,3),(2,3),(2,4),(3,4),(4,5),(5,3)\}$.
Let $K=\{1,2,3\}$ be the set of types of flow, $\tilde
U^{1}=\{(1,2),(1,3),(2,3)\}, \tilde U^{2}=\tilde
U^{3}=\{(2,3),(2,4),(3,4),(4,5),(5,3)\}$ - the sets of arcs
carrying the flow of type $k, k \in K$. We construct the networks
$S^{k}=(I^{k},U^{k}), k \in K$ (Figure \ref{graph}).

\begin{equation*}
\begin{array}{c}
x^{1}_{12}+x^{1}_{13}=4
\\
x^{1}_{23}-x^{1}_{12}=6
\\
-x^{1}_{13}-x^{1}_{23}=-10
\\
\end{array}
\end{equation*}

\begin{equation*}
\begin{array}{c}
x^{2}_{23}+x^{2}_{24}=5
\\
x^{2}_{34}-x^{2}_{23}-x^{2}_{53}=-5
\\
x^{2}_{45}-x^{2}_{24}-x^{2}_{34}=1
\\
x^{2}_{53}-x^{2}_{45}=-1
\\
\end{array}
\tag{1a}
\end{equation*}

\begin{equation*}
\begin{array}{c}
x^{3}_{23}+x^{3}_{24}=5
\\
x^{3}_{34}-x^{3}_{23}-x^{3}_{53}=-7
\\
x^{3}_{45}-x^{3}_{24}-x^{3}_{34}=1
\\
x^{3}_{53}-x^{3}_{45}=1
\\
\end{array}
\end{equation*}

\begin{equation*}
\begin{array}{c}
2x^{1}_{12}+3x^{1}_{13}+x^{1}_{23}+4x^{2}_{23}+2x^{3}_{23}
+3x^{2}_{24}-4x^{3}_{24}+2x^{2}_{34}+x^{3}_{34}-\\
x^{2}_{45}+7x^{3}_{45} +x^{2}_{53}+2x^{3}_{53}=69\\
x^{1}_{12}+2x^{1}_{13}+2x^{1}_{23}+5x^{2}_{23}+3x^{3}_{23}
-x^{2}_{24}-x^{3}_{24}+x^{2}_{34}+x^{3}_{34}-\\
2x^{2}_{45}+3x^{3}_{45}+2x^{2}_{53}-x^{3}_{53}=58\\
\end{array}
\tag{2a}
\end{equation*}

\begin{equation*}
\begin{array}{c}
x^{2}_{24}+x^{3}_{24}=1 \\
\end{array}
\tag{3a}
\end{equation*}

We choose a support of the network $S=(I,U)$  for the system (1a).
By Theorem 2 (Network Support Criterion), we build spanning trees
\\*
 $U^{k}_{T},k\in K=\{1,2,3\} $:
 $U^{1}_{T}=\{(1,2)^{1},(1,3)^{1}\}$, $U^{2}_{T}=\{(2,3)^{2},(2,4)^{2},(4,5)^{2}\}$,
$U^{3}_{T}=\{(2,4)^{3},(3,4)^{3},(4,5)^{3}\}$.

Now, we compute the set
$\{\delta^{k}(\tau,\rho),(\tau,\rho)^{k}\in U^{k}\setminus
U^{k}_{T}\}$  of characteristic vectors with respect to the
constructed spanning tree  $U^{k}_{T},k\in K=\{1,2,3\}$.
\medskip

Table 1 The set of characteristic vectors with respect to the
spanning tree $U^{1}_{T}$

\begin{tabular}{|c|c|c|c|}  \hline
$(i,j)^{1}$ &$ (1,2)^{1}$ & $(1,3)^{1}$ & $(2,3)^{1}$
\\[1mm] \hline
$\delta^{k}_{ij}(\tau,\rho)=\delta^{1}_{ij}(2,3)$ & $1 $ & $-1  $
& $1$
\\[1mm]
\hline
\end{tabular}
\\[1mm]

\medskip
Table 2 The set of characteristic vectors with respect to the
spanning tree $U^{2}_{T}$

\begin{tabular}{|c|c|c|c|c|c|}  \hline
$(i,j)^{2}$ &$ (2,3)^{2}$ & $(2,4)^{2}$ & $(4,5)^{2}$ &$
(3,4)^{2}$ &$ (5,3)^{2}$
\\[1mm] \hline
$\delta^{k}_{ij}(\tau,\rho)=\delta^{2}_{ij}(3,4)$ & $1 $ & $-1  $
& $0$ &$1$ &$0$
\\[1mm]
\hline $\delta^{k}_{ij}(\tau,\rho)=\delta^{2}_{ij}(5,3)$ & $-1 $ &
$1 $ & $1$ & $0 $& $1 $
\\[1mm]
\hline
\end{tabular}
\\[1mm]
\medskip

Table 3 The set of characteristic vectors with respect to the
spanning tree $U^{3}_{T}$

\begin{tabular}{|c|c|c|c|c|c|}  \hline
$(i,j)^{3}$ &$ (2,4)^{3}$ & $(3,4)^{3}$ & $(4,5)^{3}$ &$
(2,3)^{3}$ &$ (5,3)^{3}$
\\[1mm] \hline
$\delta^{k}_{ij}(\tau,\rho)=\delta^{3}_{ij}(2,3)$ & $-1 $ & $1 $ &
$0$ &$1$ &$0$
\\[1mm]
\hline $\delta^{k}_{ij}(\tau,\rho)=\delta^{3}_{ij}(5,3)$ & $0 $ &
$1 $ & $1$ & $0 $& $1 $
\\[1mm]
\hline
\end{tabular}
\\[1mm]

Let's compute the partial solution of the system (1a) for each
$k\in K=\{1,2,3\}$ according to the Remark 2:
$\tilde{x}^{1}=(\tilde{x}^{1}_{12},\tilde{x}^{1}_{13},\tilde{x}^{1}_{23})^{T}=
(-6,10,0)^{T}$,
$\tilde{x}^{2}=(\tilde{x}^{2}_{23},\tilde{x}^{2}_{24},\tilde{x}^{2}_{45},\tilde{x}^{2}_{34},
\tilde{x}^{2}_{53})^{T}=(5,0,1,0,0)^{T}$,
$\tilde{x}^{3}=(\tilde{x}^{3}_{24},\tilde{x}^{3}_{34},\tilde{x}^{3}_{45},
\tilde{x}^{3}_{23},\tilde{x}^{3}_{53})^{T}=(5,-7,-1,0,0)^{T}.$

We form the set $U_{C}=\displaystyle \bigcup^{3}_{k=1}
U^{k}_{C}=\{(2,3)^{1},(3,4)^{2},(2,3)^{3}\}$ of cyclic arcs. The
remaining arcs will be included into the set $U_{N}=\displaystyle
\bigcup^{3}_{k=1} U^{k}_{N}=\{(5,3)^{2},(5,3)^{3}\}$. Structures,
representing the union of the sets $U^{k}_{T}\bigcup U^{k}_{C}, k
\in K=\{1,2,3\}$ are shown on Figure 2.

\begin{figure}[h]
    \centering
        \includegraphics[width=100pt,height=100pt]{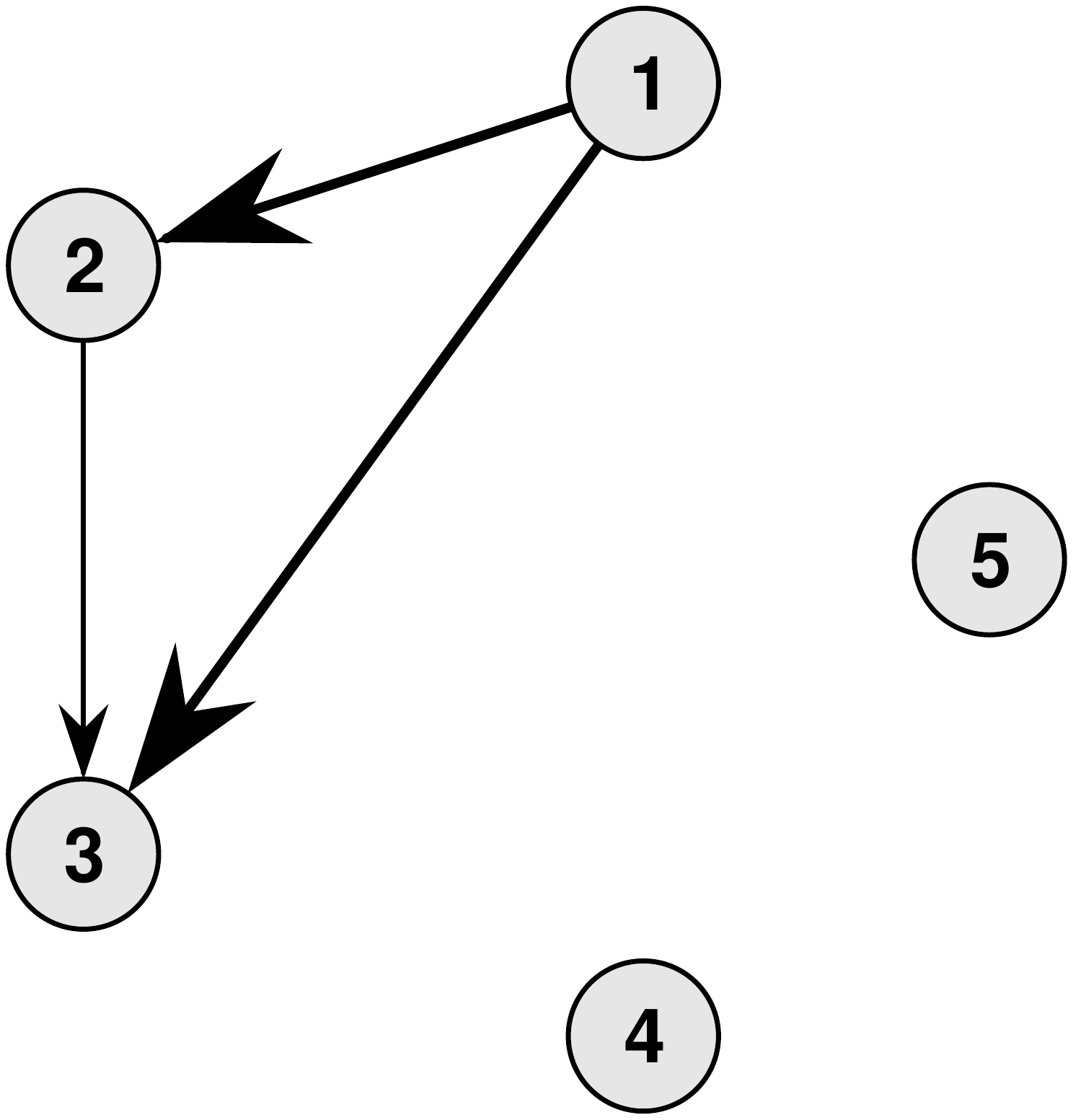}
        \includegraphics[width=100pt,height=100pt]{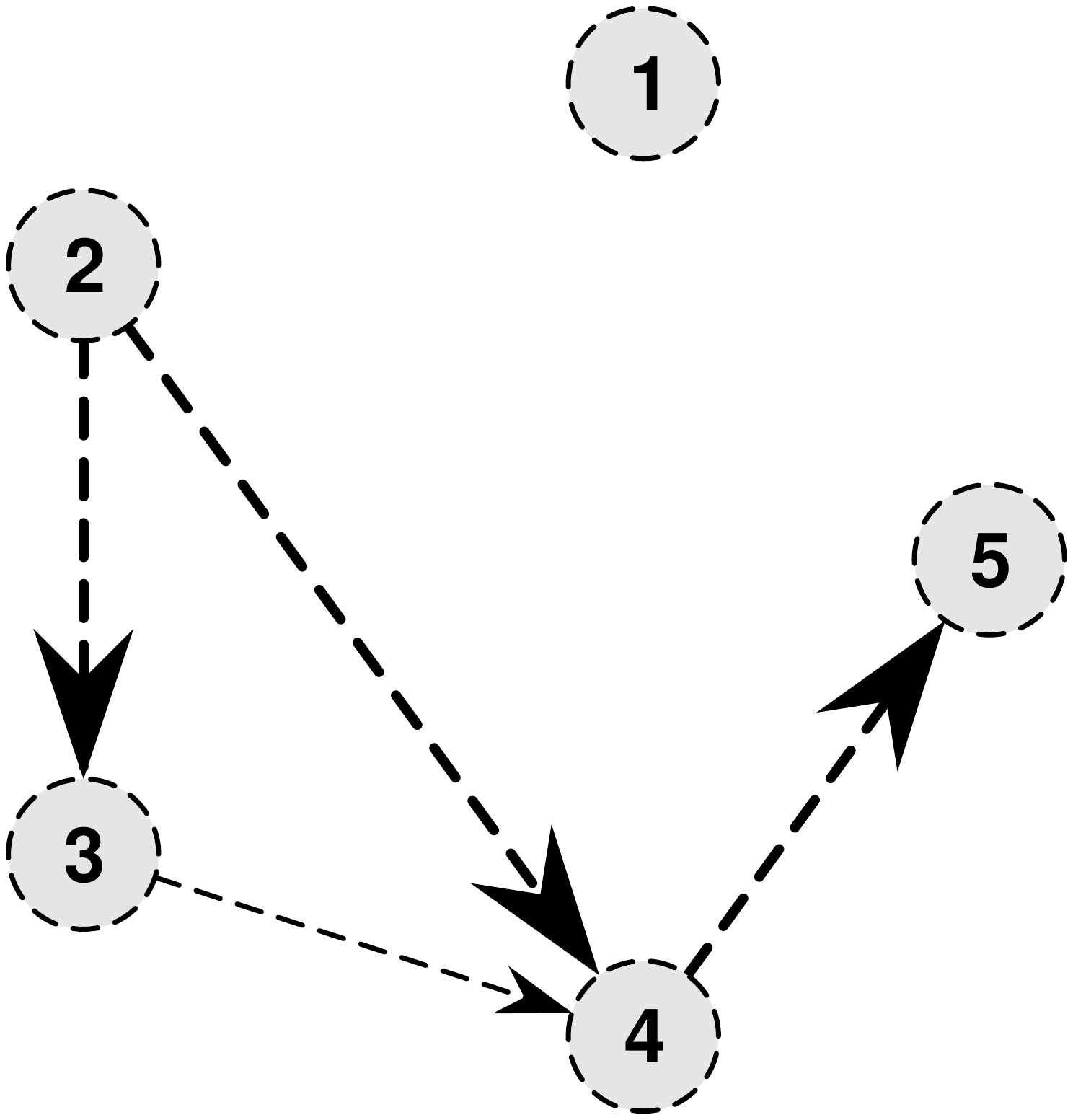}
        \includegraphics[width=100pt,height=100pt]{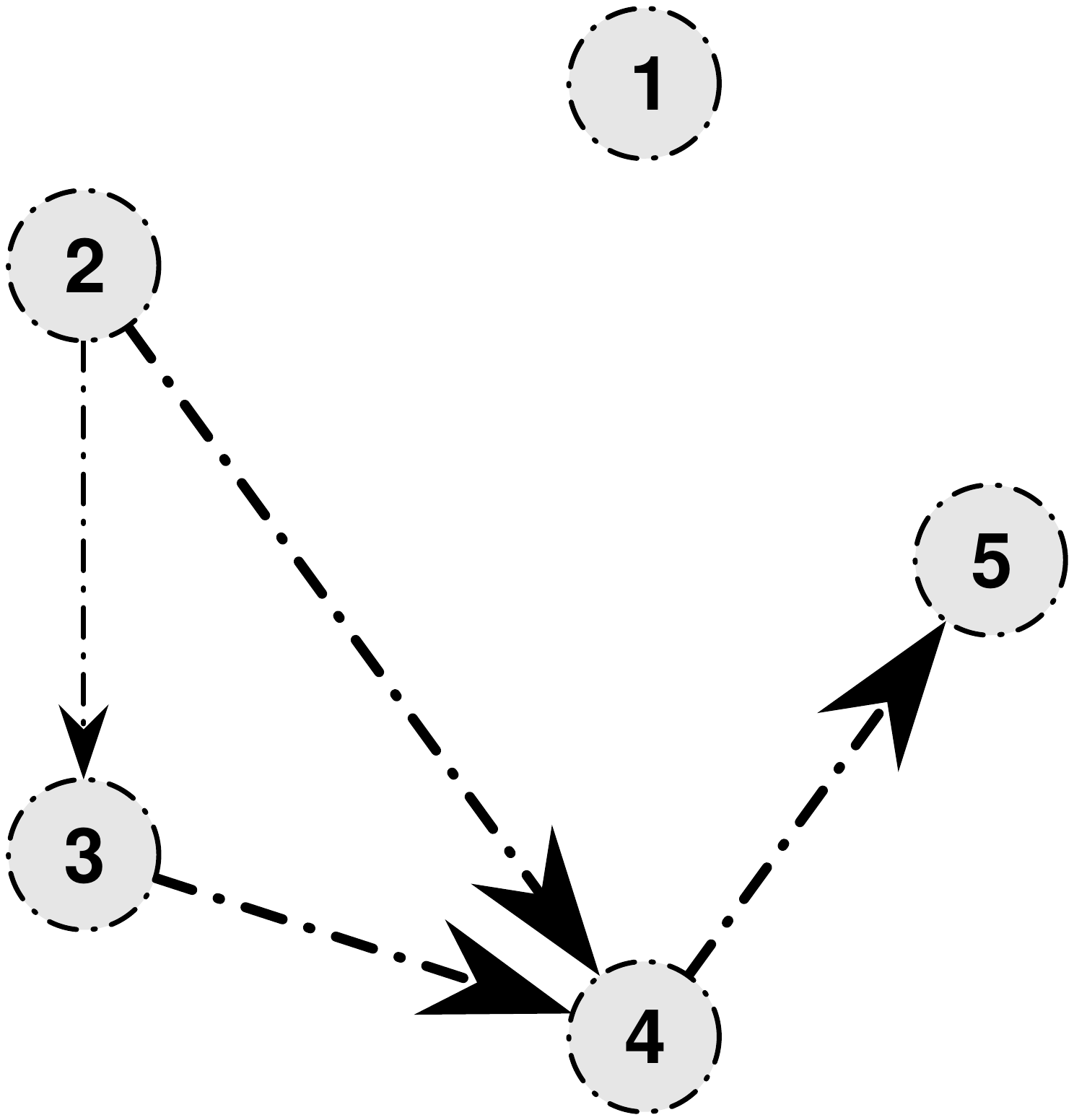}
    \caption{Sets $U^{k}_{T}\bigcup U^{k}_{C}$ for networks
    $S^{k}, k \in K=\{1,2,3\}$}
    \label{trees}
\end{figure}

Using formula (\ref{eq18}) we compute the determinants of the
cycles  $L^{k}_{\tau\rho}$, entailed by the arcs
$(\tau,\rho)^{k}\in U^{k}\backslash U^{k}_{T}$, for each  $k\in
K=\{1,2,3\}$, with respect to the equation (2a) with the number
$p=1,2$ (Table 4).
\medskip

Table 4  Determinants of the cycles $L^{k}_{\tau\rho}$, entailed
by the arcs $(\tau,\rho)^{k}\in U^{k}\backslash U^{k}_{T}, k\in
K=\{1,2,3\}$

\begin{tabular}{|c|c|c|c|c|c|}  \hline
$(\tau,\rho)^{k}$ &$ (2,3)^{1}$ & $(3,4)^{2}$ & $(5,3)^{2}$ &$
(2,3)^{3}$ &$ (5,3)^{3}$
\\[1mm] \hline
$R_{1}(L^{k}_{\tau\rho})$ & $0 $ & $3 $ & $-1$ &$7 $ &$10$
\\[1mm]
\hline $R_{2}(L^{k}_{\tau\rho})$ & $1 $ & $7 $ & $-6 $ & $5 $& $3
$
\\[1mm]
\hline
\end{tabular}
\\[1mm]

Now, let's compute the values $\delta_{ij}(L^{k}_{\tau\rho}),
(i,j)\in U_{0}, (\tau,\rho)^{k} \in U^{k}\setminus U^{k}_{T}, k\in
K=\{1,2,3\}$ according to the formula (\ref{eq24}) for the example
(1a)-(3a), $U_{0}=\{(2,4)\}$, $K_{0}(2,4)=\{2,3\}$ (Table 5).
\medskip

Table 5  The values  $\delta_{ij}(L^{k}_{\tau\rho}), (i,j)\in
U_{0}, (\tau,\rho)^{k} \in U^{k}\setminus U^{k}_{T}, k\in
K=\{1,2,3\}$

\begin{tabular}{|c|c|c|c|c|c|}  \hline
$(\tau,\rho)^{k}$ &$ (2,3)^{1}$ & $(3,4)^{2}$ & $(5,3)^{2}$ &$
(2,3)^{3}$ &$ (5,3)^{3}$
\\[1mm] \hline
$\delta_{24}(L^{k}_{\tau\rho})$ & $0 $ & $-1 $ & $1$ &$-1 $ &$0$
\\[1mm]
\hline
\end{tabular}
\\[1mm]

Before assembling the matrix $D$ of the system (\ref{eq28}), let's
number the arcs of the set  : $U_{C}=
\{(2,3)^{1},(3,4)^{2},(2,3)^{3}\}:t(2,3)^{1}=1,\,t(3,4)^{2}=2,\,
t(2,3)^{3}=3$. The numbering within the set $U_{0}=\{(2,4)\}$  is
trivial: $\xi(2,4)=1$.

First, we construct the matrix
$D_{1}=(R_{p}(L^{k}_{\tau\rho}),\,p=\overline{1,2},\,t(\tau,\rho)^{k}=\overline{1,3})$
of the determinants of the cycles $L^{k}_{\tau\rho}$ , entailed by
the arcs $(\tau,\rho)^{k}\in U_{C}$, by selecting the
corresponding columns from the Table 4:
\begin{displaymath}
D_{1}=\left(\begin{array}{ccc}
0&3&7 \\
1&7&5
\end{array} \right).
\end{displaymath}

Similarly, by selecting the corresponding columns from the Table
5, we form the matrix
$D_{2}=(\delta_{24}(L^{k}_{\tau\rho}),\,\xi(2,4)=1,\,
t(\tau,\rho)^{k}=\overline{1,3})$:
$$D_{2}=(0\quad -1\quad -1).$$

Thus, joining $D_{1}$  and $D_{2}$  together, we obtain the matrix
of the system (\ref{eq28}):

\begin{displaymath}
D=\left(\begin{array}{ccc}
0&3&7 \\
1&7&5 \\
0&-1&-1
\end{array} \right),\mbox{det} D\neq 0.
\end{displaymath}

Let us compute the vector $\beta$  in the right hand side of
(\ref{eq28}) using formulae (\ref{eq29}):

$$\beta_{1}=A^{1}-R_{1}(L^{2}_{53})x^{2}_{53}-R_{1}(L^{3}_{53})x^{3}_{53},$$
$$\beta_{2}=A^{2}-R_{2}(L^{2}_{53})x^{2}_{53}-R_{2}(L^{3}_{53})x^{3}_{53},$$
$$\beta_{3}=A_{24}-\delta_{24}(L^{2}_{53})x^{2}_{53}-\delta_{24}(L^{3}_{53})x^{3}_{53}.$$

The values  $R_{p}(L^{2}_{53})$,  $R_{p}(L^{3}_{53})$, $p=1,2$  of
the determinants of the cycles  $L^{k}_{\tau\rho}$, entailed by
the arcs  $(\tau,\rho)^{k}\in U_{N}$, as well as the values
$\delta_{24}(L^{2}_{53})$ and $\delta_{24}(L^{3}_{53})$, are
already computed and stored within the Table 4 and Table 5. The
numbers $A^{1},A^{2},A_{24}$ are evaluated using the formulae
(\ref{eq19}) and (\ref{eq26}):
$$A^{1}=\alpha_{1}-\lambda_{12}^{11}\tilde{x}^{1}_{12}
-\lambda_{13}^{11}\tilde{x}^{1}_{13}
-\lambda_{23}^{21}\tilde{x}^{2}_{23}
-\lambda_{24}^{21}\tilde{x}^{2}_{24}
-\lambda_{45}^{21}\tilde{x}^{2}_{45}
-\lambda_{24}^{31}\tilde{x}^{3}_{24}
-\lambda_{34}^{31}\tilde{x}^{3}_{34}
-\lambda_{45}^{31}\tilde{x}^{3}_{45} =66,$$
$$A^{2}=\alpha_{2}-\lambda_{12}^{12}\tilde{x}^{1}_{12}
-\lambda_{13}^{12}\tilde{x}^{1}_{13}
-\lambda_{23}^{22}\tilde{x}^{2}_{23}
-\lambda_{24}^{22}\tilde{x}^{2}_{24}
-\lambda_{45}^{22}\tilde{x}^{2}_{45}
-\lambda_{24}^{32}\tilde{x}^{3}_{24}
-\lambda_{34}^{32}\tilde{x}^{3}_{34}
-\lambda_{45}^{32}\tilde{x}^{3}_{45} =36,$$
$$A_{24}=z_{24}-\tilde{x}^{2}_{24}-\tilde{x}^{3}_{24}=-4.$$

Thus, we have defined the vector
\begin{math}
\beta=\left(\begin{array}{c}
66+x^{2}_{53}-10x^{3}_{53}\\
36+6x^{2}_{53}-3x^{3}_{53}\\
-4-x^{2}_{53}
\end{array}
\right).
\end{math}

Since the matrix $D$ turned out to be non-singular, we can use
formula (\ref{eq30}) for finding the solution
$x_{C}=(x^{k}_{\tau\rho},(\tau,\rho)^{k}\in U_{C}^{k}, k \in K)$
of the system (\ref{eq28}):
\begin{displaymath}
\left(\begin{array}{c}
x^{1}_{23} \\
x^{2}_{34} \\
x^{3}_{23} \\
\end{array} \right)=
\left( \begin{array}{rcr}
\displaystyle \frac{1}{2}&1&\displaystyle \frac{17}{2} \\
-\displaystyle \frac{1}{4}&0&-\displaystyle \frac{7}{4}\\
\displaystyle \frac{1}{4}&0&\displaystyle \frac{3}{4}\\
\end{array}
\right)
\left(\begin{array}{c}
66+x^{2}_{53}-10x^{3}_{53}\\
36+6x^{2}_{53}-3x^{3}_{53}\\
-4-x^{2}_{53}\\
\end{array}
\right).
\end{displaymath}

Finally, using formulae (\ref{eq31}) - (\ref{eq32}), we can define
the solution of the system (1a)-(3a) with $x^{2}_{53},x^{3}_{53}$,
being independent variables:
\begin{displaymath}
x^{1}_{23}=35-2x^{2}_{53}-8x^{3}_{53},x^{2}_{34}=-\frac{19}{2}+
\frac{3}{2}x^{2}_{53}+\frac{5}{2}x^{3}_{53},
x^{3}_{23}=\frac{27}{2}-
\frac{1}{2}x^{2}_{53}-\frac{5}{2}x^{3}_{53},
\end{displaymath}

$$
x^{1}_{12}=29-2x^{2}_{53}-8x^{3}_{53},
$$
$$
x^{1}_{13}=-25+2x^{2}_{53}+8x^{3}_{53},
$$

$$
x^{2}_{23}=-\frac{9}{2}+
\frac{1}{2}x^{2}_{53}+\frac{5}{2}x^{3}_{53},
$$
$$
x^{2}_{24}=\frac{19}{2}-
\frac{1}{2}x^{2}_{53}-\frac{5}{2}x^{3}_{53},
$$
$$
x^{2}_{45}=x^{2}_{53}+1,
$$

$$
x^{3}_{24}=-\frac{17}{2}+
\frac{1}{2}x^{2}_{53}+\frac{5}{2}x^{3}_{53},
$$
$$
x^{3}_{34}=\frac{13}{2}-
\frac{1}{2}x^{2}_{53}-\frac{3}{2}x^{3}_{53},
$$
$$
x^{3}_{45}=x^{3}_{53}-1,
$$
$$
x^{2}_{53},x^{3}_{53}\in \mathbf{R}.
$$

\bigskip
{\footnotesize \rm
\begin{tabular}{ll}
L.A. Pilipchuk           &\qquad E.S. Vecharynski\\
Belarus State University &\qquad University of Colorado Denver\\
Department of Applied Mathematics 
 &\qquad Department of Mathematical Sciences\\
 and Computer Science
&\qquad Campus Box 170\\
Nezalezhnasci Ave. 4 &\qquad P.O. Box 173364\\
 Minsk         &\qquad Denver, CO 80217-3364 USA\\    
220050 Belarus      &\qquad e-mail: \\
e-mail: pilipchuk@bsu.by 
      &\qquad Yaugen.Vecharynski@ucdenver.edu\\
\end{tabular}
}
\end{document}